\documentclass[12pt]{article}
\usepackage{latexsym}
\usepackage{amssymb}
\usepackage{graphics}
\newtheorem{prop}{Proposition}[section]
\newtheorem{lem}[prop]{Lemma}

\newcommand{\bmin}[1]{\hspace{-1mm}\stackrel{\frown}{\hspace{2mm}#1}\hspace{-1mm}}
\newcommand{\bmax}[1]{\hspace{-1mm}\stackrel{\frown}{#1\hspace{2mm}}\hspace{-1mm}}
\newcommand{\adj}[2]{\stackrel{\frown}{#1 #2}}
\newcommand{\bmaxmin}[1]{\stackrel{\frown\frown}{\hspace{.75mm}#1}}
\newcommand{\badj}[2]{\stackrel{\hspace{-.3cm}\resizebox{1cm}{2mm}{$\frown$}}{#1
#2}}
\newcommand{\badji}[2]{\stackrel{\hspace{-2mm}\resizebox{4mm}{1.5mm}{$\frown$}}{#1
#2}}
\begin{document}

\bibliographystyle{elsart-num-sort}
\title{Avoidance of Partitions of a Three-element Set}
\author{Adam M. Goyt\footnote{This work was partially done while the author was visiting DIMACS.}\\Department of Mathematics\\ Michigan State University\\ East Lansing, Michigan 48824-1027\\ goytadam@msu.edu\\www.math.msu.edu/$\sim$goytadam}
\date{May 5, 2006}

\maketitle

\begin{abstract}

Klazar defined and studied a notion of pattern avoidance for set
partitions, which is an analogue of pattern avoidance for
permutations. Sagan considered partitions which avoid a single
partition of three elements.  We enumerate partitions which avoid
any family of partitions of a 3-element set as was done by Simion
and Schmidt for permutations. We also consider even and odd set
partitions.  We provide enumerative results for set partitions
restricted by generalized partition patterns, which are an analogue
of the generalized permutation patterns of Babson and
Steingr{\'{\i}}msson. Finally, in the spirit of work done by Babson
and Steingr{\'{\i}}msson, we will show how these generalized
partition patterns can be used to describe set partition statistics.

\end{abstract}


\medskip


\section{Introduction}

Pattern avoidance in permutations was first introduced by Knuth in
\cite{Knuthvol3}, and is currently an area of very active research.
An approach to studying pattern avoidance and containment that deals
with set partitions was introduced and studied by Klazar in
\cite{abbafree, KlazarPartI, KlazarPartII} and continued by Sagan
in~\cite{Saganpartitionpatterns}. The extensively studied non-crossing
partitions defined by Kreweras \cite{KrewerasNonCross} can be viewed
as those which avoid a specific pattern with four elements. We will
focus on the enumeration of all partitions of an $n$-element set
which avoid a family of partitions of a 3-element set. To make these
notions of pattern containment for permutations and set partitions
precise and to see their connections we will need some definitions.

If $f:S\rightarrow T$ is a function from set $S$ to set $T$, then
$f$ acts element-wise on objects constructed from $S$. For example,
if $a_1a_2\dots a_n$ is a permutation of elements of $S$ then
$f(a_1a_2\dots a_n)=f(a_1)f(a_2)\dots f(a_n)$.  Also, define $[n]$
to be the set $\{1,2,\dots,n\}$ and $[k,n]$ to be the set
$\{k,k+1,\dots n\}$.

Suppose that $S\subseteq\mathbb{Z}$ is a set with $\#S=n$,  then the
{\it standardization} map corresponding to $S$ is the unique order
preserving bijection $St_S:S\rightarrow [n]$.  For example if
$S=\{2,5,7,10\}$ then $St_S(2)=1$, $St_S(5)=2$, $St_S(7)=3$, and
$St_S(10)=4$.  When it is clear from context what set the
standardization map is acting on, we will omit the subscript $S$.

Let $p=a_1a_2\dots a_k\in S_k$ be a given permutation, called the
{\it pattern}, where $S_k$ is the symmetric group on $k$ letters. A
permutation $q=b_1b_2\dots b_n\in S_n$ {\it contains} the pattern
$p$ if there is a subsequence $ q'=b_{i_1}b_{i_2}\dots b_{i_k}$ of
$q$ with $St(q')=p$. Otherwise $q$ {\it avoids} $p$. For example the
permutation $q=32145$ contains 6 copies of the pattern 213, namely
324, 325, 314, 315, 214, and 215. On the other hand $q$ avoids the
pattern 132.  For $R\subseteq S_k$, let

$$S_n(R)=\{q\in S_n:\mbox{$q$ avoids every pattern }p\in R\}.$$ The
problem of enumerating $S_n(R)$  for $R\subseteq S_3$ was considered
by Simion and Schmidt~\cite{SimionSchmidt}.  We will consider the
analogous problem for patterns in partitions.

A {\it partition} $\pi$ of set $S\subseteq\mathbb{Z}$, written
$\pi\vdash S$, is a family of nonempty, pairwise disjoint subsets
$B_1,B_2,\dots,B_k$ of $S$ called {\it blocks} such that
$\bigcup_{i=1}^k B_i=S$.  We write $\pi=B_1/B_2/\dots/B_k$ and
define the {\it length} of $\pi$, written $\ell(\pi)$, to be the number
of blocks. Since the order of the blocks does not matter, we will
always write our partitions in the {\it canonical order} where
$$\min B_1 < \min B_2 <\dots<\min B_k.$$  We will also always write
the elements of each block in increasing order.  For example,
$137/26/45\vdash[7]$ has length 3.

Let $$\Pi_n=\{\pi\vdash [n]\}$$ be the set of all partitions of
$[n]$. Suppose $\sigma$ is a set partition of length $m$ and $\pi$
is a partition of length $\ell$. Then $\sigma$ {\it contains} $\pi$,
written $\pi\subseteq \sigma$, if there are $\ell$ different blocks of
$\sigma$ each containing a block of $\pi$. For example
$\sigma=137/26/45$ contains $\pi =2/37/5$ but does not contain
$\pi'=2/37/6$ because 2 and 6 are in the same block of $\sigma$.

Let $\pi\in\Pi_k$ be a given set partition called the {\it pattern}.
 A partition $\sigma\in\Pi_n$ {\it contains the pattern} $\pi$ if
there is some $\sigma'\subseteq\sigma$ with $St(\sigma')=\pi$.
Otherwise $\pi$ {\it avoids} $\sigma$.  For example
$\sigma=137/26/45$ contains six copies of the pattern $\pi=14/2/3$,
namely 17/2/4, 17/2/5, 17/4/6, 17/5/6, 26/3/4, and 26/3/5. It is
important to note here that when looking for a copy of $\pi$ in
$\sigma$, the order of the blocks does not matter.  On the other
hand consider the pattern $\pi'=1/234$.  To be contained in $\sigma$
the copy of the block 234 of $\pi'$ must be contained in a block of
size three or larger. The only such block of $\sigma$ is 137.  It is
impossible to find an element smaller than 1, so $\sigma$ does not
contain a copy of $\pi'$.  For $R\subseteq\Pi_k$, let
$$\Pi_n(R)=\{\sigma\in\Pi_n:\mbox{$\sigma$ avoids every pattern }\pi\in
R\}.$$

The set of non-crossing partitions mentioned above may be defined as
the set $\Pi_n(13/24)$. It is known that $\#\Pi_n(13/24)=C_n$, where
$C_n$ is the $n$th Catalan number \cite{KrewerasNonCross},
\cite{RStanleyVol2}. For a survey of results about non-crossing
partitions see Simion's paper \cite{SimionSurvey}.

Sagan~\cite{Saganpartitionpatterns} has provided enumerative results for
$\Pi_n(R)$ when $\#R=1$. In the spirit of work done by Simion and
Schmidt on permutation patterns~\cite{SimionSchmidt}, we will
enumerate $\Pi_n(R)$ for $\#R\geq 2$.  We then define the sign of a
partition and enumerate the set of signed partitions of $[n]$
avoiding particular patterns. In section 5, we define generalized
patterns analogous to the generalized permutation patterns of Babson
and Steingr{\'{\i}}msson~\cite{BabsonSteingrimsson}, and provide
enumerative results for those.  Finally, we will show how these
generalized partition patterns can be used to describe set partition
statistics.

\section{Double Restrictions}

In this section we will consider the case of $\#\Pi_n(R)$ where
$\#R=2$.  Given a set partition $\sigma=B_1/B_2/\dots/B_k\vdash
[n]$, let $\sigma^c=B_1^c/B_2^c/\dots/B_k^c$ be the {\it complement}
of $\sigma$ where $$B_i^c=\{n-a+1:a\in B_i\}.$$
 For example if $\sigma = 126/3/45$ then $\sigma^c=156/23/4$. The
following result is obvious, so we omit the proof.

\begin{prop}[Sagan]
For $n\geq1$,
\begin{eqnarray*}
\Pi_n(\sigma^c)&=&\{\pi^c:\pi\in \Pi_n(\sigma)\},\\
\#\Pi_n(\sigma^c)&=&\#\Pi_n(\sigma).\: 
\end{eqnarray*}

\end{prop}
The following Lemma is an immediate consequence of Proposition 2.1.

\begin{lem}

\begin{eqnarray*}
\#\Pi_n(12/3, 123)& =& \#\Pi_n(1/23, 123)\\
\#\Pi_n(1/2/3, 12/3)& =& \#\Pi_n(1/2/3, 1/23)\\
\#\Pi_n(12/3, 13/2)&= &\#\Pi_n(1/23, 13/2).\: 
\end{eqnarray*}
\end{lem}

There are 10 different sets $R$ with elements from $\Pi_3$ and
$\#R=2$, so by Lemma 1 there are seven different cases to consider.
 Note that $\#\Pi_0=1$ by letting the empty set partition itself.
Since any partition in $\Pi_1$ or $\Pi_2$ cannot possibly contain a
partition of [3], we have $\#\Pi_0(R)=1$, $\#\Pi_1(R)=1$ and
$\#\Pi_2(R)=2$ for all $R\subseteq \Pi_3$.  The fact that
$\#\Pi_3=5$ implies that $\#\Pi_3(R)=3$ for any $R\subset \Pi_3$,
with $\#R=2$. Hence, it suffices to consider $n\geq4$ in the
following results.

A partition $\sigma\vdash [n]$ is {\it layered} if $\sigma$ is of
the form $[1,i]/[i+1,j]/[j+1,k]/\dots/[\ell+1,n]$.  An example of a
layered partition is $\sigma=123/4/56/789$.  A partition $\sigma$ is
a {\it matching} if $\#B\leq2$ for every block $B$ of $\sigma$.

We will use the following results of Sagan~\cite{Saganpartitionpatterns}
repeatedly, so we state them now.

\begin{prop}[Sagan]
\begin{eqnarray}
\label{1/2/3}\Pi_n(1/2/3)&=&\{\sigma:l(\sigma)\leq2\},\\
\Pi_n(12/3)&=&\{\sigma=B_1/B_2/\dots /B_k : {\min}B_i=i \:for\:
each\:
i,\: and\nonumber\\
\label{12/3}&&\hspace{4.2cm}
[k+1,n]\subseteq B_i \:for \:some \:i\},\\
\label{13/2}\Pi_n(13/2)&=&\{\sigma:\sigma\: is\: layered\},\\
\label{123}\Pi_n(123)&=&\{\sigma:\sigma\: is\: a\: matching\}.\:
\square
\end{eqnarray}
\end{prop}

\begin{prop}  For all $n\geq3$,
\begin{eqnarray*}
\Pi_n(1/2/3, 12/3)& =& \{12\dots n,\: 1/23\dots n, \:13\dots
n/2\},\\
\# \Pi_n(1/2/3, 12/3)& =& 3.\\
\end{eqnarray*}
\end{prop}

\noindent{\bf Proof}:  Let $\sigma \in \Pi_n(1/2/3, 12/3)$. By
(\ref{1/2/3}), $\sigma$ may have at most two blocks.  If $\ell(\sigma)
= 1$ then $\sigma = 12\dots n$. If $\ell(\sigma) = 2$ then by
(\ref{12/3}), we must have
$[3,n]\subset B_i$ for $i=$ 1 or 2. $\square$\\

\begin{prop}  For all $n\geq1$,
\begin{eqnarray*}\Pi_n(1/2/3,13/2)& = &\{\sigma : \sigma=12\dots k/ (k+1)(k+2)\dots n\: for\: some\: k \in
[n]\},\\ \#\Pi_n(1/2/3,13/2)&=&n.\\
\end{eqnarray*}
\end{prop}

\noindent {\bf Proof}:  If $\sigma \in \Pi_n(1/2/3,13/2)$ then
$\sigma$ is layered by (\ref{13/2}), and $\ell(\sigma)\leq2$ by
(\ref{1/2/3}). Hence $\sigma$ is of the form described above.  The
enumeration follows immediately. $\square$

\begin{prop}  \begin{eqnarray*}
\Pi_n(1/2/3, 123) &=&\left\{\begin{array}{cr}
\{12/34,13/24,14/23\} & n=4,\\
\emptyset & n\geq 5.
\end{array} \right.\\
\#\Pi_n(1/2/3, 123)& =& \left\{
\begin{array}{lr}
\hspace{1.9cm}3 & \hspace{1.8cm}n=4,\\
\hspace{1.9cm}0 & \hspace{1.8cm}n\geq5.
\end{array}\right.\\
\end{eqnarray*}
\end{prop}

\noindent{\bf Proof}:  If $n\geq5$ and $\sigma \vdash[n]$, then
$\ell(\sigma)\geq3$ or $\sigma$ has a block of size $\geq3$ by the
Pigeonhole Principle. Thus by (\ref{1/2/3}) and (\ref{123}),
$\Pi_n(1/2/3,123)=\emptyset$ for $n\geq5$. The case $n=4$ is easy to
check.  $\square$

\begin{prop}  For all $n\geq3$,\\
\begin{eqnarray*}
\Pi_n(1/23,12/3) &=& \{12\dots n,\: 1/2/\dots/ n,\:
1n/2/3/\dots/n-1\}, \\
\#\Pi_n(1/23,12/3) &=& 3.\\
\end{eqnarray*}
\end{prop}

\noindent {\bf Proof:}  Let $\sigma=B_1/B_2/\dots /B_k$ avoid 12/3.
If $k=1$ then $\sigma=12\dots n$, which avoids 1/23.  Similarly,
when $k=n$, we have $\sigma=1/2/\dots /n$, which avoids 1/23. If
$k=n-1$ and $n\in B_i$ for $i\geq2$ then $B_1/B_i$ is a copy of
$1/23$. Thus $n\in B_1$ and $\sigma=1n/2/3/\dots/n-1$.  If $1<k<
n-1$ then, by (\ref{12/3}), we must have $\{n-1,n\} \subseteq B_i$
for some $i$, and there is at least one more block. Hence $\sigma$
contains a copy of $1/23$, and so this case can not occur.
$\square$

\begin{prop}  For all $n\geq1$,
\begin{eqnarray*}\Pi_n(12/3, 13/2)
&= &\{\sigma = 1/2/\dots/k-1/k(k+1)\dots n,\: for \: some \: k\in[n]\},\\
\#\Pi_n(12/3, 13/2) &=& n.\\
\end{eqnarray*}
\end{prop}

\noindent{\bf Proof}:  Suppose $\sigma=B_1/B_2/\dots/
B_k\in\Pi_n(12/3,13/2)$.  Then by (\ref{12/3}) we have $i\in B_i$
for each $i$ and exactly one of the $B_i$ contains $[k+1,n]$.  From
(\ref{13/2}) we have that $\sigma$ must be layered.  So $[k+1,n]\in
B_k$, and $B_k=[k,n]$.  Thus there is exactly one
$\sigma\in\Pi_n(12/3,13/2)$ of length $k$ for each $k\in[n]$.
$\square$

\begin{prop}  For all $n\geq1$,
\begin{eqnarray*}
\Pi_n(12/3, 123)& =& \{\sigma = B_1/B_2/\dots /B_k : \min B_i=i,\:
and\:
k=n-1\: or \: n\},\\
\# \Pi_n(12/3, 123)& =& n.\\
\end{eqnarray*}
\end{prop}

\noindent{\bf Proof}:  Assume $\sigma=B_1/B_2/\dots/B_k \in
\Pi_n(12/3, 123)$.  Then by (\ref{12/3}) and (\ref{123}), $k=n-1$ or
$n$. The result follows. $\square$

\medskip

Let $F_n$ be the $n^{th}$ Fibonacci number, initialized by $F_0=1$
and $F_1=1$.  A {\it composition} of an integer $n$ is an ordered
collection of positive integers $n_1,n_2,\dots ,n_k$ such that
$n=n_1+n_2+\dots +n_k$. The $n_i$ are called {\it parts}.  It is
easy to see that $F_n$ counts the number of compositions of $n$ with
parts of size 1 or 2.

\begin{prop} For all $n\geq0$,
\begin{eqnarray*}\Pi_n(13/2, 123)&=&\{\sigma : \sigma\mathrm{\:is\:a\:layered\:matching}\},\\
\#\Pi_n(13/2, 123)&=&F_{n}.\\ \end{eqnarray*}

\end{prop}

\noindent{\bf Proof}:  Any $\sigma\in\Pi_n(13/2,123)$ must be
layered by (\ref{13/2}) and a matching by (\ref{123}).

There is a bijection between the compositions of $n$ with parts of
size 1 or 2 and the partitions of $[n]$ that are layered matchings.
If $\sigma\in\Pi_n(13/2,123)$ and $\sigma=B_1/B_2/\dots /B_k$, then
we map $\sigma$ to the composition $n=n_1+n_2+\dots+n_k$ with
$n_i=\#B_i$.  $\square$

\medskip

From the results above we know that
$$\#\Pi_n(1/2/3,13/2)=\#\Pi_n(12/3, 13/2)=\#\Pi_n(12/3, 123)=n,$$ and
we have a very nice description of the elements in each of these
sets.  It is interesting to note that one gets similar results when
avoiding certain sets of permutations in $S_3$.

\begin{prop}[Simion, Schmidt] For every $n\geq1$,
\begin{eqnarray*}
\#S_n(123,132,231)&=&\#S_n(123,213,312)=n.\\
\#S_n(132,231,321)&=&\#S_n(213,312,321)=n.\\
\end{eqnarray*}
And:
\begin{eqnarray*}
q\in S_n(123,132,231)&\iff& q=(n,n-1,\dots,k+1,k-1,k-2,\dots,2,1,k),\\
q\in S_n(123,213,312)&\iff& q=(n,n-1,\dots,k+1,1,2,3,\dots,k),\\
q\in S_n(132,231,321)&\iff& q=(n-1,n-2,\dots k+1,n,k,k-1,\dots,2,1),\\
q\in S_n(213,312,321)&\iff& q=(k-1,\dots,3,2,1,n,n-1,\dots,k).\:\square\\
\end{eqnarray*}
\end{prop}

The Fibonacci numbers also occur when avoiding permutations.

\begin{prop}[Simion, Schmidt]  For every $n\geq 1$,
$$\#S_n(123,132,213)=F_{n}.\: \square$$
\end{prop}

There is a simple map $\Phi:\Pi_n\rightarrow S_n$, given by sending
$\sigma=B_1/B_2/\dots/B_k$ to $B_kB_{k-1}\dots B_1$. For example,
$\Phi(1/23/4/56)= 564231$.

\begin{prop} The map $\Phi$ restricts to a bijection from the set $\Pi_n(13/2,123)$ to the set $S_n(123,132,213)$.
\end{prop}

\noindent{\bf Proof:}  We may describe $q\in S_n(123,132,213)$
recursively.  To avoid the patterns 123 and 213, we must have
$q^{-1}(n)\leq 2$.  If $q^{-1}(n)=1$ then the remaining positions
form a permutation in $S_{n-1}(123,132,213)$.  If $q^{-1}(n)=2$ then
$q^{-1}(n-1)=1$, otherwise there will be a copy of 132 in $q$.  The
remaining positions form a permutation in $S_{n-2}(123,132,213)$.

Suppose $\sigma=B_1/B_2/\dots/B_k\in\Pi_n(13/2,123)$, then
$B_k=\{n\}$ or $\{n-1,n\}$.  The permutation $\Phi(\sigma)$ thus
begins with $n$ or $n-1,n$.  Inductively, one can see that this
restriction of the map $\Phi$ is well defined.

To prove that the restricted $\Phi$ is a bijection we provide its
inverse map. Let $q=q_1q_2\dots q_n\in S_n(123,132,213)$ then we say
that $q_k$ is a {\it descent} if $q_k>q_{k+1}$.  Let
$D=\{q_{i_1},q_{i_2},\dots, q_{i_{\ell}}\}$ be the set of descents of
$q$, with $i_1<i_2<\dots <i_{\ell}$.  Then
$$\Phi^{-1}(q)=q_{i_{\ell}+1}q_{i_{\ell}+2}\dots q_n/q_{i_{\ell-1}+1}\dots
q_{i_{\ell}}/\dots/q_1\dots q_{i_1}.$$

For example $\Phi^{-1}(564231)=1/23/4/56$ because its descent set is
$D=\{3,4,6\}.$

We now show that $\Phi^{-1}$ is well defined.  Every $q\in
S_n(123,132,213)$ must have a descent in at least one of its first
two positions.  After this initial descent there may be no more than
one position between any two descents.  Thus the blocks of
$\Phi^{-1}(q)$ will have size at most 2, and from the description of
the elements of $S_n(123,132,213)$ above $\Phi^{-1}(q)$ will be
layered.

The fact that $\Phi$ and $\Phi^{-1}$ are inverses follows easily
from the descriptions of the maps. $\square$

\section{Higher Order Restrictions}

We begin, as with double restrictions, by reducing the number of
cases.  The following Lemma is a consequence of Proposition 2.1.

\begin{lem} 
\begin{eqnarray*}
\#\Pi_n(1/2/3, 12/3,123) &=& \#\Pi_n(1/2/3, 1/23,123),\\
\#\Pi_n(1/2/3, 12/3,13/2)&=&\#\Pi_n(1/2/3, 1/23,13/2),\\
\#\Pi_n(12/3,13/2,123)&=&\#\Pi_n(1/23,13/2,123).\:\square
\end{eqnarray*}
\end{lem}

The results for $\#\Pi_n(R)$ where $\#R=3$ are easy to prove. Table
3.3 describes these sets and gives their enumeration for $n\geq4$.
The following proposition describes $\#\Pi_n(R)$ for $\#R\geq4$.  We
omit the simple proof.

\begin{prop}  For $R\subseteq\Pi_3$ with $\#R\geq4$ and $n\geq4$,
$$\#\Pi_n(R)=\left\{
\begin{array}{ll}
0&if\:\{1/2/3,123\}\subseteq R,\\
1& else. \: \square
\end{array}
\right.
$$
\end{prop}

\bigskip

\begin{tabular}{|c|c|c|}
\hline $R$ & $\Pi_n(R)$ & $\#\Pi_n(R)$\\
\hline $\{1/2/3,12/3,13/2\}$& $\{12\dots n,1/23\dots n\}$&2\\
\hline $\{1/2/3,12/3,123\}$&$\emptyset$&0\\
\hline $\{1/2/3,13/2, 123\}$&$\{12/34\}$ & 1 if $n=4$\\
&$\emptyset$&0 if $n\geq5$\\
\hline $\{1/2/3,1/23,12/3\}$&$\{12\dots n\}$&1\\
\hline
$\{12/3,13/2,123\}$&$\{1/2/\dots/n,1/2/\dots/n-2/(n-1)n\}$&2\\
\hline $\{1/23,12/3,13/2\}$&$\{123\dots n,1/2/\dots/n\}$&2\\
\hline $\{1/23,12/3,123\}$&$\{1/2/\dots/n,1n/2/3/\dots/n-1\}$&2\\
\hline
\end{tabular}

\bigskip

 {\bf Table 3.3:} Enumeration of partitions restricted by 3 patterns

\section{Even and Odd Set Partitions}

In this section we will consider the number of even and odd
partitions of the set $[n]$, which avoid a single pattern of length
three.  A partition $\sigma\vdash [n]$ with $l(\sigma)=k$ has {\it
sign}, $$\mathrm{sgn}(\sigma)=(-1)^{n-k}.$$  {\it Even} partitions
$\sigma$ satisfy $\mathrm{sgn}(\sigma)=1$, and {\it odd} partitions
$\sigma$ satisfy $\mathrm{sgn}(\sigma)=-1$.
We will use the following notation:\\
\begin{eqnarray*}
E\Pi_n(\pi)&=&\{\sigma \vdash [n]: \mathrm{sgn}(\sigma)=1\},\\
O\Pi_n(\pi)&=&\{\sigma \vdash [n]: \mathrm{sgn}(\sigma)=-1\}.\\
\end{eqnarray*}

The following follows directly from the definitions.

\begin{lem}  The sign of $\sigma$ is the same as the
sign of $\sigma^c$.  Thus $\#E\Pi_n(12/3)=\#E\Pi_n(1/23)$ and
$\#O\Pi_n(12/3)=\#O\Pi_n(1/23).$ $\square$
\end{lem}

We will use the following result of Sagan \cite{Saganpartitionpatterns}
repeatedly, so we state it now.  Define the {\it double factorial}
by $$(2i)!!=1\cdot3\cdot5\cdots(2i-1).$$

\begin{prop}[Sagan] 
\begin{eqnarray}
\label{enum1/2/3}\#\Pi_n(1/2/3)&=&2^{n-1},\\
\label{enum12/3}\#\Pi_n(12/3)&=&{n \choose 2} +1,\\
\label{enum13/2}\#\Pi_n(13/2)&=&2^{n-1}.\\
\label{enum123}\#\Pi_n(123)&=&\sum_{i=0}^{\lfloor n/2
\rfloor}{n\choose 2i}(2i)!!\:\square
\end{eqnarray}
\end{prop}

We now consider single restrictions.  By Lemma 4.1 there are only
four cases.

\begin{prop} For all odd $n\geq1$,\\
\begin{eqnarray*}
\#E\Pi_n(1/2/3)&=& 1,\\
\#O\Pi_n(1/2/3)&=& 2^{n-1}-1.\\
\end{eqnarray*}
For all even $n\geq2$,
\begin{eqnarray*}
\#E\Pi_n(1/2/3)&=& 2^{n-1}-1,\\\
\#O\Pi_n(1/2/3)&=& 1.\\
\end{eqnarray*}
\end{prop}

\noindent{\bf Proof:}  By (\ref{1/2/3}), any $\sigma\in\Pi_n(1/2/3)$
must have $\ell(\sigma)\leq2$.  If $n$ is odd then a partition of
length 1 will be even and a partition of length 2 will be odd. There
is only one partition of length 1, and
$\#O\Pi_n(\pi)+\#E\Pi_n(\pi)=\#\Pi_n(\pi)$ for any pattern $\pi$.
Thus, the result holds for odd $n$ by (\ref{enum1/2/3}). The proof
for even $n$ is similar. $\square$

\begin{prop} For all odd $n\geq0$,\\
\begin{eqnarray*}
\#E\Pi_n(12/3)&=&\left\lfloor\frac{(n-1)^2}{4}\right\rfloor+1,\\
\#O\Pi_n(12/3)&=&\left\lfloor\frac{n^2}{4}\right\rfloor.
\end{eqnarray*}
\end{prop}

\noindent{\bf Proof}:  By (\ref{12/3}) we have, for $n$ odd,
$$\#E\Pi_n(12/3)=1+\sum_{k=0}^{\frac{n-3}{2}}(2k+1)=1+\frac{(n-1)^2}{4}=\left\lfloor\frac{(n-1)^2}{4}\right\rfloor+1,$$\\
and by (\ref{enum12/3})
$$\#O\Pi_n(12/3)={n\choose2}+1-\frac{(n-1)^2}{4}-1=\left\lfloor\frac{n^2}{4}\right\rfloor.$$
The proof for even $n$ is similar.  $\square$

\begin{prop} For all $n\geq1$,
$$\#O\Pi_n(13/2)=\#E\Pi_n(13/2)=2^{n-2}.$$
\end{prop}

\noindent{\bf Proof}:  By (\ref{enum13/2}) it suffices to give a
sign reversing involution $\psi:\Pi_n(13/2)\rightarrow\Pi_n(13/2)$.
By (3), $\sigma\in\Pi_n(13/2)$ is layered, so it is of the form
$\sigma=B_1/B_2/\dots/B_k$, where either $B_k=\{n\}$ or
$B_k\supset\{n\}$.  Let
\begin{eqnarray*}
\psi(\sigma)&=&\left\{\begin{array}{lr}
B_1/B_2/\dots/B_{k-1}\cup\{n\}&\mathrm{if}\:B_k=\{n\},\\
B_1/B_2/\dots/B_k-\{n\}/n&\mathrm{if}\: B_k\supset\{n\}.
\end{array}\right.
\end{eqnarray*}
Notice that $\psi(\sigma)$ is still layered for any
$\sigma\in\Pi_n(13/2)$, so $\psi$ is well defined.  And, $\psi$ is
its own inverse because it either moves $n$ into the block preceding
it if $\{n\}$ is a block and into its own block otherwise. Also,
$\psi$ changes the sign of $\sigma$ by either increasing or
decreasing the length of $\sigma$ by 1. $\square$

\begin{prop} For all $n\geq1$, \\
\begin{eqnarray*}
\#E\Pi_n(123)&=&\sum_{i=0}^{\lfloor\frac{n-2}{4}\rfloor}{n\choose{4i+2}}(4i+2)!!\\
\#O\Pi_n(123)&=&\sum_{i=0}^{\lfloor\frac{n}{4}\rfloor}{n\choose{4i}}(4i)!!\\
\end{eqnarray*}
\end{prop}

\noindent{\bf Proof}:  Any $\sigma \in\Pi_n(123)$ is a matching. If
$i$ blocks of $\sigma$ have 2 elements each and the remaining blocks
are singletons then $\sigma$ has $i+(n-2i)= n-i$ blocks. Thus
$\mathrm{sgn}(\sigma)=(-1)^{n-(n-i)}=(-1)^{i}$.  So the even and odd
counts are obtained by taking the appropriate terms from
(\ref{enum123}). $\square$

\medskip

Table 4.7 gives the results for $\#E\Pi_n(R)$ and $\#O\Pi_n(R)$
where $\#R\geq2$ and $n\geq4$.  We prove the enumeration of
$E\Pi_n(13/2,123)$ and $O\Pi_n(13/2,123)$ as an example and leave
the rest to the reader.

\begin{prop}  

\begin{eqnarray*}
\#E\Pi_n(13/2,123)&=&\left\{\begin{array}{lr}\lceil F_n/2 \rceil &
\mathrm{for} \:n\equiv0,1 \:(\mathrm{mod} \: 6),\\
F_n/2 & \mathrm{for} \: n \equiv2,5 \:(\mathrm{mod} \: 6),\\
\lfloor F_n/2 \rfloor & \mathrm{for} \: n\equiv 3,4 \:(\mathrm{mod} \:
6).\end{array}\right.\\
\#O\Pi_n(13/2,123)&=&\left\{\begin{array}{lr}\lfloor F_n/2 \rfloor &
\mathrm{for} \:n\equiv0,1 \:(\mathrm{mod} \: 6),\\
F_n/2 & \mathrm{for} \:n\equiv2,5 \:(\mathrm{mod} \: 6),\\
\lceil F_n/2 \rceil & \mathrm{for} \: n\equiv 3,4 \:(\mathrm{mod} \:
6).\end{array}\right.
\end{eqnarray*}
\end{prop}

\bigskip

\begin{tabular}{|c|c|c|}
\hline $R$ & $\#E\Pi_n(R)$ & $\#O\Pi_n(R)$ \\ \hline
$\{1/2/3,12/3\}$ & 1 for $n$ odd & 2 for $n$ odd \\
 & 2 for $n$ even & 1 for $n$ even \\ \hline
$\{1/2/3, 13/2\}$ & $1$ for $n$ odd & $n-1$ for $n$ odd \\
 & $n-1$ for $n$ even & 1 for $n$ even \\ \hline
 $\{1/2/3, 123\}$ & 3 for $n=4$ & 0 \\
 &  0 for $n\geq5$ &\\ \hline
 $\{1/23,12/3\}$ & 2 for $n$ odd & 1 for $n$ odd \\
  & 1 for $n$ even & 2 for $n$ even\\ \hline
 $\{12/3,13/2\}$ & $\lceil n/2\rceil$ & $\lfloor n/2 \rfloor$  \\ \hline
 $\{12/3,123\}$ & 1 & $n-1$ \\ \hline
 $\{13/2, 123\}$  & $\lceil F_n/2 \rceil$ for $n\equiv 0,1 $ $(\mathrm{mod}\: 6)$ & $\lfloor F_n/2 \rfloor$ for $n\equiv 0,1$ $ (\mathrm{mod}\: 6)$\\
 & $F_n/2$ for $n\equiv 2,5$ $(\mathrm{mod}\: 6)$ & $F_n/2$ for $n\equiv 2,5$ $ (\mathrm{mod}\: 6)$\\
  & $\lfloor F_n/2 \rfloor$ for $n\equiv 3,4 $ $(\mathrm{mod}\: 6)$ & $\lceil F_n/2 \rceil$ for $n\equiv 3,4 $ $(\mathrm{mod}\: 6)$\\ \hline
 $\{1/2/3, 1/23, 12/3\}$  & 1 for $n$ odd & 0 for $n$ odd\\
  & 0 for $n$ even & 1 for $n$ even\\ \hline
 $\{1/2/3, 12/3, 13/2\}$ &1& 1 \\ \hline
 $\{1/2/3, 12/3, 123\}$ & 0 & 0 \\ \hline
 $\{1/2/3, 13/2, 123\}$ & 1 & 0 \\ \hline
 $\{1/23, 12/3, 13/2\}$ & 2 for $n$ odd & 0 for $n$ odd \\
 & 1 for $n$ even & 1 for $n$ even \\ \hline
 $\{1/23, 12/3, 123\}$ & 1 & 1 \\ \hline
 $\{12/3, 13/2, 123\}$ & 1 & 1 \\ \hline
 $\{1/2/3, 1/23, 12/3, 13/2\}$ & 1 for $n$ even & 0 for $n$ odd \\
 & 0 for $n$ odd & 1 for $n$ even \\ \hline
$\{1/2/3, 1/23, 12/3, 123\}$ & 0 & 0 \\ \hline $\{1/2/3, 12/3, 13/2,
123\}$ & 0 & 0 \\ \hline
$\{1/23, 12/3, 13/2, 123\}$ & 1 & 0 \\
\hline $\{1/2/3, 1/23, 12/3, 13/2, 123\}$ & 0 & 0 \\ \hline

\end{tabular}

\bigskip

{\bf Table 4.7:} Enumeration of even and odd partitions restricted
by at least 2 patterns

\bigskip

\noindent{\bf Proof:}  Let $\sigma=B_1/B_2/\dots/B_k\in
\Pi_n(13/2,123)$. Then $B_k=\{n\}$ or $\{n-1,n\}$.  If $B_k=\{n\}$
then $B_1/B_2/\dots/B_{k-1}$ is a layered matching of $[n-1]$ and
$\mathrm{sgn}(B_1/B_2/\dots/B_{k-1})=\mathrm{sgn}(\sigma)$. If
$B_k=\{n-1,n\}$ then $B_1/B_2/\dots/B_{k-1}$ is a layered matching
of $[n-2]$ and
$\mathrm{sgn}(B_1/B_2/\dots/B_{k-1})=-\mathrm{sgn}(\sigma)$. Thus we
have that
$$\#E\Pi_n(13/2,123)=\#E\Pi_{n-1}(13/2,123)+\#O\Pi_{n-2}(13/2,123).$$
Similarly,
$$\#O\Pi_n(13/2,123)=\#O\Pi_{n-1}(13/2,123)+\#E\Pi_{n-2}(13/2,123).$$

Now induct on $n$.  To show that the proposition is true when $0\leq
n\leq 5$ is easy. This leaves us with twelve cases to check for the
inductive step. We will show one of them. It is easy to see that
$F_n$ is odd unless $n\equiv 2,5\: (\mathrm{mod}\:6)$.

Suppose that $n\equiv 4 \:(\mathrm{mod}\:6)$.  Then we have
\begin{eqnarray*}
\#E\Pi_n(13/2,123)&=&
\#E\Pi_{n-1}(13/2,123)+\#O\Pi_{n-2}(13/2,123)\\&=&\lfloor
F_{n-1}/2 \rfloor+F_{n-2}/2\\
&=&\frac{F_{n-1}-1+F_{n-2}}{2}\\
&=&\lfloor F_n/2 \rfloor.\:\square
\end{eqnarray*}

\section{Generalized Partition Patterns}

Babson and Steingr{\'{\i}}msson~\cite{BabsonSteingrimsson} defined
generalized patterns for permutations.  These were patterns in which
certain elements were required to be consecutive.  Generalized
permutation patterns were used to describe permutation statistics
and classify Mahonian statistics.   In this section we will define a
similar notion for set partition patterns and consider the avoidance
case. In the next section we will show that generalized partition
patterns can be used to describe set partition statistics.

Recall that if $\sigma=B_1/B_2/\dots/B_k$ is a partition then the
blocks are written in such a way that $\min B_1<\min B_2< \dots
<\min B_k$.  This gives us a well defined notion of adjacency of
blocks, where we consider $B_i$ as being adjacent to both $B_{i-1}$
and $B_{i+1}$. Consider the partition $\sigma= 147/25/36$ and the
pattern $\pi=13/2$. Suppose now that a copy of $\pi$ must appear in
adjacent blocks. Then 17/2 is still a copy, but 17/3 is not. We may
also have the blocks in the restricted copy of 13/2 in the opposite
order making 25/4 a copy of $\pi$ in $\sigma$. We will denote $\pi$
with the adjacency restriction by the {\it generalized pattern}
$\rho=13|2$. In general, we will denote block adjacency using a
vertical bar.

Recall that the elements of a block are put in order by size, which
gives us a way to consider adjacent elements.  Now, suppose we want 
to find a copy of 13/2 in $\sigma=147/25/36$, but we require that
the elements that represent 1 and 3 in this copy are adjacent. In
this case 14/3 is a copy of 13/2, but 17/6 is not, since 1 and 7 are
not adjacent in their block.  We will denote this by the {\it
generalized pattern} $\rho=\adj{1}{3}\hspace{-1mm}/2$. In general,
we will denote element adjacency by placing an arc over the
elements, which must be adjacent.

If $\rho$ is a generalized pattern, then the notation $\Pi_n(\rho)$
denotes the set of partitions of $[n]$, which avoid $\rho$.
Similarly, if $R$ is any set of generalized patterns then $\Pi_n(R)$
is the set of partitions of $[n]$, which avoid all generalized
patterns in $R$.

We are interested in enumerating the $\Pi_n(R)$ where $R$ is a set
of partitions of $[3]$ at least one of which contains an adjacency
restriction. It turns out that the adjacency restrictions do not
actually restrict most of the original patterns. This is summed up
in the next lemma.

\begin{lem}  The following are true for generalized patterns:

$\begin{array}{rcccccl}
\Pi_n(1/2/3)&=&\Pi_n(1|2/3)&=&\Pi_n(1/2|3)&=&\Pi_n(1|2|3),\\
\label{b}\Pi_n(1/23)&=&\Pi_n(1|23)&=&\Pi_n(1/\hspace{-1mm}\adj{2}{3})&=&\Pi_n(1|\hspace{-1mm}\adj{2}{3}),\\
\label{c}\Pi_n(13/2)&=&\Pi_n(\adj{1}{3}\hspace{-1mm}/2)&=&\Pi_n(13|2)&=&\Pi_n(\adj{1}{3}\hspace{-1mm}|2),\\
\label{d}\Pi_n(123)&=&\Pi_n(\adj{1}{2}\hspace{-1mm}3)&=&\Pi_n(1\hspace{-1mm}\adj{2}{3})&=&\Pi_n(\adj{1}{2}\hspace{-1mm}\stackrel{\frown}{\hspace{1mm}{3}}),\\
\label{e}\Pi_n(12/3)&=&\Pi_n(\adj{1}{2}\hspace{-1mm}/3),&&&&\\
\label{f}\Pi_n(12|3)&=&\Pi_n(\adj{1}{2}\hspace{-1mm}|3).&&&&
\end{array}$

\end{lem}

\noindent{\bf Proof:}  We will only prove the second line as the
others are very similar.  First we show that
$\Pi_n(1/23)=\Pi_n(1|23)$.  It is obvious that if a partition
$\sigma\vdash[n]$ contains a copy of $1|23$ then it contains a copy
of 1/23.  So it will suffice to show the other containment holds.
Let $\sigma=B_1/B_2/\dots/B_k\vdash[n]$ contain a copy $a/bc$ of
1/23. Suppose $a\in B_s$ and $b,c\in B_t$. If $s<t$ then the block
$B_{t-1}$ exists and $\min B_{t-1}<\min B_t\leq b<c$. Letting
$d=\min B_{t-1}$ gives a copy $d/bc$ of $1|23$ in $\sigma$.  If
$s>t$ then $B_{t+1}$ exists and $\min B_{t+1}\leq a<b<c$.  Letting
$e=\min B_{t+1}$ gives a copy $e/bc$ of $1|23$ in $\sigma$.  We
remind the reader that the adjacent blocks of the copy of $1|23$ may
appear in either order in $\sigma$.

Now we will show that
$\Pi_n(1/23)=\Pi_n(1/\hspace{-1mm}\adj{2}{3})$.  Again, it suffices
to show that if $\sigma\vdash[n]$ contains a copy of 1/23 then it
contains a copy of $1/\hspace{-1mm}\adj{2}{3}$.  Given a copy $a/bc$
of 1/23 in $\sigma$, if $b$ and $c$ are not adjacent in their block
$B$ then let $d$ be the minimum of all of the elements of $B$ which
are larger than $b$. Thus $a/bd$ is a copy of
$1/\hspace{-1mm}\adj{2}{3}$ in $\sigma$.  These two observations can
be used to prove the remaining equality. $\square$

\medskip

Let $R$ be a set of generalized patterns, and let $S$ be the same
set with adjacency restrictions dropped.  That is if, for example,
$1|\hspace{-1mm}\adj{2}{3}\hspace{1mm}\in R$ then $1/23\in S$, and
$S$ only contains patterns without adjacency restrictions. Lemma 5.1
says that unless $12|3$ or $\adj{1}{2}\hspace{-1mm}|3\in R$, we have
that $\Pi_n(R)=\Pi_n(S)$. However, since we have
$\Pi_n(12|3)=\Pi_n(\adj{1}{2}|3)$,  we only need to consider cases
when $12|3\in R$. The sets $\Pi_n(S)$ were enumerated in sections 2
and 3, so we need only enumerate the sets $\Pi_n(S\cup\{12|3\})$
where $S\subseteq\Pi_3-\{12/3\}$.

\begin{prop} Let $S\subseteq\Pi_3-\{12/3\}$ then
$\Pi_n(S\cup\{12|3\})=\Pi_n(S\cup\{12/3\})$ unless $S=\emptyset$ or
$\{123\}$.

\end{prop}

\noindent{\bf Proof:}  The cases where $\#S\geq2$ follow
automatically from those with $\#S=1$ and Lemma 5.1. The three cases
with $\#S=1$ are very similar, so we will only prove the statement
for $S=\{13/2\}$.  Let $\sigma\in\Pi_n(13/2,12|3)$, then $\sigma$
must be layered. Thus any copy of $12/3$ in $\sigma$ easily reduces
to a copy of $12|3$ as in the proof of Lemma 5.1. $\square$

\medskip

The following lemma describes the elements of $\Pi_n(12|3)$.

\begin{lem}  We have $\sigma\in \Pi_n(12|3)$ if and only if whenever a block $B_t$ of
$\sigma$ satisfies $\#B_t\geq2$, then
$$\#B_{t-1}=1\:and\:\#B_{t+1}=1.$$  Furthermore, if
$B_{t+1}=\{a\}$ then $a<b$ for every $b\in B_{t}-\{\min B_t\}$.

\end{lem}

\noindent {\bf Proof:}  First we show that
$\sigma=B_1/B_2/\dots/B_k\in\Pi_n(12|3)$ can be described as above.
Let $\#B_t\geq2$ and suppose that $B_{t-1}$ contains at least 2
elements and let $a<b$ be the two smallest elements of $B_{t-1}$.
Let $c<d$ be the two smallest elements of $B_t$. By the definition
of canonical order, $a<c$. If $b<d$, then $ab/d$ is a copy of
$12|3$. If $b>d$, then $cd/b$ is a copy of $12|3$ another
contradiction. The proof that $\#B_{t+1}=1$ is similar. The single
element in $B_{t+1}$ must be larger than $c$ by definition.  If it
is larger than any other element of $B_t$ we will again have an
unwanted copy of $12|3$.

Now, suppose that $\sigma\in\Pi_n$ has the structure described
above.  Then it is straight forward to show that $\sigma$ cannot contain a copy of $12|3$.  $\square$

\medskip

First we will consider the case where $S=\emptyset$ in Proposition
5.2.  Let $a_n=\#\Pi_n(12|3)$ and let
$$f(x)=\sum_{n\geq0}a_n\frac{x^n}{n!}$$ be the corresponding exponential
generating function.

\begin{prop}  For $n\geq2$, $$a_{n}=a_{n-1}+1+\sum_{k=1}^{n-2}{n-2\choose k}a_{n-k-2}$$ with the
initial conditions $a_0=1$ and $a_1=1$, and $f(x)$ satisfies the
differential equation $$y''=y'+y(e^x-1)+e^x.$$
\end{prop}

\noindent{\bf Proof:}  That $\#\Pi_0(12|3)=\#\Pi_1(12|3)=1$ is
obvious. Let $\sigma=B_1/B_2/\dots /B_k \in\Pi_n(12|3)$. Either
$\#B_1=1$ or $\#B_1\geq2$.  If $\#B_1=1$ then, by the definition of
canonical order, $B_1=\{1\}$. Clearly any $12|3$ avoiding partition
of the set $[2,n]$ will still avoid $12|3$ if we prepend the block
$\{1\}$. This gives the first term of the recursion.

Now suppose that $\#B_1\geq2$, then either $\sigma=12\dots n$ or
not. The case where $\sigma=12\dots n$ is counted by the 1 in the
recursion. If $\sigma\not=12\dots n$ then, by Lemma 5.3, we must
have $B_2=\{2\}$. If $k$ of the elements from $[3,n]$ are in $B_1$,
then the remaining $n-k-2$ elements must form a $12|3$ avoiding
partition. This establishes the recursion.

Using the recursion to produce the differential equation satisfied
by $f(x)$ is routine and is left the reader. $\square$

\medskip

The substitution $y=ue^{x/2}$ simplifies the equation to
$$u''=u(e^x-\frac{3}{4})+e^{x/2}.$$

Using Maple, we obtain the solution
$$u=C_1\cdot I_{\sqrt{-3}}(2e^{x/2})+C_2\cdot
K_{\sqrt{-3}}(e^{x/2})+$$ $$2I_{\sqrt{-3}}(2e^{x/2})\int
K_{\sqrt{-3}}(e^{x/2}e^{x/2})dx -$$ $$ 2K_{\sqrt{-3}}(e^{x/2})\int
I_{\sqrt{-3}}(2e^{x/2})e^{x/2}dx,$$ for certain constants $C_1$ and
$C_2$, where $I_n(z)$ and $K_n(z)$ are the modified Bessel functions
of the first and second kinds respectively. There are known
combinatorial interpretations for certain Bessel functions. See, for
example, \cite{Bessel1} and \cite{Bessel2}. It is unlikely, however,
that there is a combinatorial interpretation for the Bessel
functions appearing in the exponential generating function
$f(x)=ue^{-x/2}$, since $K_{\sqrt{-3}}(e^{x/2})$ is not well defined
as a formal power series.

Now, we turn our focus to $\Pi_n(123,12|3)$.  Let $b_n=\#\Pi_n(123,
12|3)$ and $$g(x)=\sum_{n\geq0}b_n\frac{x^n}{n!}$$ be the
corresponding exponential generating function.

The proof of the following Proposition is very similar to the proof of Proposition 5.4 and is omitted.

\begin{prop}  For $n\geq3$, $$b_n=b_{n-1}+(n-2)b_{n-3}$$ with the initial conditions $b_0=1$, $b_1=1$, and
$b_2=2$.  Also, $g(x)$ satisfies the differential equation
$$y'''=y''+xy'+y.\: \square$$
\end{prop}

Using Maple, we obtain the solution
$$y=D_1e^{x/2}Ai(1/4+x+D_2e^{x/2}Bi(1/4+x)+$$
$$\hspace{1.6cm}D_3e^{x/2}\left(Ai(1/4+x)\int Bi(1/4+x)e^{-x/2}dx-\right.$$
$$\hspace{3.2cm}\left.\int Ai(1/4+x)e^{-x/2}dx
Bi(1/4+x)\right),$$ for constants $D_1$, $D_2$, and $D_3$, where
$Ai$ and $Bi$ are Airy functions.

It is not terribly surprising that Airy functions appear, since
these functions are closely related to Bessel functions and
$\Pi_n(123,12|3)$ is a subset of the set $\Pi_n(12|3)$. There do not
seem to be any existing combinatorial interpretations of Airy
functions. There is also unlikely to be a combinatorial
interpretation of this generating function due to the fact that
$Ai(1/4+x)$ is not well defined as a formal power series.

For completeness we will consider the cases where odd and even set
partition avoid generalized set partitions.  As before only the
cases $O\Pi_n(R)$ and $E\Pi_n(R)$ where $R=\{12|3\}$ or
$\{123,12|3\}$ are new.

Let $oa_n=\#O\Pi_n(12|3)$ and $ea_n=\#E\Pi_n(12|3)$.  Let
$ob_n=\#O\Pi_n(123,12|3)$ and $eb_n=\#E\Pi_n(123,12|3)$.  The
following propositions easily follow from the recursions above.  We
let $\chi$ be the truth function, where $\chi$ of a statement is 1
if the statement is true and 0 if the statement is false.

\begin{prop}  For $n\geq2$,
$$oa_n=oa_{n-1}+\chi(n\:is\:even)+\sum_{l=2,\hspace{1mm} l \hspace{1mm} even}^{n-2}{n-2\choose l}oa_{n-2-l}+\sum_{l=1,\hspace{1mm} l \hspace{1mm} odd}^{n-2}{n-2 \choose l}ea_{n-2-l},$$
and
$$ea_n=ea_{n-1}+\chi(n\:is\:odd)+\sum_{l=2,\hspace{1mm} l \hspace{1mm} even}^{n-2}{n-2\choose l}ea_{n-2-l}+\sum_{l=1,\hspace{1mm} l \hspace{1mm} odd}^{n-2}{n-2 \choose l}oa_{n-2-l}.\:\square$$
\end{prop}

\begin{prop}  For $n\geq3$

$$ob_n=ob_{n-1}+(n-2)eb_{n-3},$$
and
$$eb_n=eb_{n-1}+(n-2)ob_{n-3}.\:\square$$

\end{prop}

\section{Set Partition Statistics}

Carlitz~\cite{CarlitzStirl1,CarlitzStirl2} and
Gould~\cite{GouldStirl} were the first to give versions of the
$q$-Stirling numbers of the second kind.  In \cite{Milnestirling},
Milne introduces an inversion and dual inversion statistic on set
partitions, whose distributions over partitions of $[n]$ with $k$
blocks produce these two $q$-Stirling numbers of the second kind.
Later, Sagan \cite{SaganMaj} introduced the major index and dual
major index of a set partition, whose distributions produced the
same two $q$-Stirling numbers of the second kind.  At around the
same time, Wachs and White \cite{WachsWhite} investigated four
natural statistics, which they called {\it lb}, {\it ls}, {\it rb},
and {\it rs}, again producing the same two $q$-Stirling numbers of
the second kind. Other statistics of interest are the number of
crossings, nestings and alignments of a partition, see for example
\cite{3cross}, \cite{VasTab}, or \cite{KasZeng}. In this section we
will show that all of these statistics can be described in the
language of generalized partition patterns.

We will need some more notation.  Consider the pattern $\pi = 1/23$.
If we are looking for a copy of $\pi$ in $\sigma=137/26/45$, but we
want the element representing 1 in the copy to be the minimum of its
block then $1/45$ is a copy, but 3/45 is not.  We will represent
this generalized pattern by $\bmin{1}/23$. And in general, we will
denote such a generalized pattern by putting an arc over the first
element of the block, in which we want the minimum to occur. In the
same fashion, if we want the element representing 1 in a copy of
1/23 to be the maximum in its block, then we denote the pattern by
$\:\bmax{1}/23$. If we want the element representing 1 in a copy of
1/23 to be both the minimum and the maximum of its block, then we
denote the pattern by $\bmaxmin{1}/23$.

In the sequel, if we say $\rho$ is a pattern then $\rho$ may or may
not have adjacency restrictions. Let $\rho$ be a pattern and
$\sigma\in\Pi_n$. Then $\rho$ will be treated as a function from
$\Pi_n$ to the nonnegative integers by letting $\rho(\sigma)$ be the
number of copies of $\rho$ in $\sigma$.  If we have patterns
$\rho_1, \rho_2,\dots, \rho_{\ell}$ then
$$(\rho_1+\rho_2+\dots+\rho_{\ell})(\sigma)=\rho_1(\sigma)+\rho_2(\sigma)+\dots+\rho_{\ell}(\sigma).$$

We begin with the inversion statistic.  Let
$\sigma=B_1/B_2/\dots/B_k\in\Pi_n$ and $b\in B_i$.  We will say that
$(b,B_j)$ is an {\it inversion} if $b> \min B_j$ and $i<j$. Define
the {\it inversion number} of $\sigma$, written
$\mathrm{inv}(\sigma)$, to be the number of inversions in $\sigma$.

We may calculate $\mathrm{inv}(\sigma)$ by summing, over all
elements $b\in [n]$, the number of inversions of the form $(b,B_j)$.
This observation leads to the next Proposition.

\begin{prop}  For any $\sigma\in\Pi_n$, $$\mathrm{inv}(\sigma)=(\hspace{1mm}\bmin{1}3/\bmin{2}\hspace{1mm})(\sigma).$$
\end{prop}

\noindent {\bf Proof:} We will show that there is a one to one
correspondence between inversions and copies of
$\bmin{1}\hspace{-1mm}3/\bmin{2}\hspace{1mm}$.  Let
$\sigma=B_1/B_2/\dots /B_k$. Let $b\in B_i$ and $(b,B_j)$ be an
inversion.  If $a=\min B_i$ and $c=\min B_j$ then $(b,B_j)$
corresponds to the copy $ab/c$ of $\bmin{1}3/\bmin{2}\hspace{1mm}$.
Conversely, if $ab/c$ is a copy of $\bmin{1}3/\bmin{2}\hspace{1mm}$,
then $a=\min B_i$ and $c=\min B_j$ where $i<j$ since $a<c$. Also,
$b>c=\min B_j$.  Thus, the copy $ab/c$ yields the inversion $(b,
B_j)$.
 $\square$

\medskip

Let $\sigma=B_1/B_2/\dots/B_k$ be a partition.  We will say that
$(b,B_{i+1})$ is a {\it descent} of $\sigma$ if $b\in B_i$ and
$b>\min B_{i+1}$. Let $d_i$ be the number of descents of $\sigma$ in block $B_i$.
Then the major index of $\sigma$ is
$$\mathrm{maj}(\sigma)=\sum_{i=1}^{k-1}id_i=d_1+2d_2+\dots+(k-1)d_{k-1}.$$

Notice that each descent $(b,B_{i+1})$ contributes $i$ to the major index.

\begin{prop} For any $\sigma\in\Pi_n$, $$\mathrm{maj}(\sigma)=(\:\bmin{1}3|\bmin{2}\hspace{1mm}+\bmin{1}/\bmin{2}4|\bmin{3}\:)(\sigma).$$
\end{prop}

\noindent {\bf Proof:} Let $\sigma=B_1/B_2/\dots/B_k$ and $b\in
B_i$.  Let $\rho_1=\bmin{1}3|\bmin{2}\hspace{1mm}$ and
$\rho_2=\bmin{1}/\bmin{2}4|\bmin{3}\hspace{1mm}$. We will first show
that $(b,B_{i+1})$ is a descent if and only if $b$ represents the 3
in a copy of $\rho_1$, or, for $i\geq2$, the 4 in a copy of
$\rho_2$.  Then we will show that each descent $(b,B_{i+1})$
contributes $i$ to the right hand side.

Let $(b,B_{i+1})$ be a descent.  If $a=\min B_i$ and $c=\min
B_{i+1}$ then $ab/c$ is a copy of $\rho_1$ where $b$ represents the
3.  If additionally $i\geq2$ and we let $d=\min B_{j}$ where $j<i$
then $d/ab/c$ is a copy of $\rho_2$, in which $b$ represents the 4.
For the converse, let $ab/c$ be a copy of $\rho_1$, then $c=\min
B_{i+1}$ for some $i$, and $(b,B_{i+1})$ is a descent. Similarly, a
copy $d/ab/c$ of $\rho_2$ with $c=\min B_{i+1}$ for some $i\geq2$
produces the descent $(b,B_{i+1})$.

If $(b,B_{i+1})$ is a descent, then there is exactly one copy of
$\rho_2$ with $b$ representing 3, since the 1 in $\rho_1$ must be
represented by $a=\min B_i$, and the 2 must be represented by
$c=\min B_{i+1}$.  Now, if $b$ represents the 4 in a copy of
$\rho_2$ then the 2 must be represented by $a=\min B_i$, and the 3
must be represented by $c=\min B_{i+1}$.  But now the 1 may be
represented by the minimum of any block appearing before $B_i$.  So
the total contribution of the two patterns is $1+(i-1)=i$.
$\square$

\medskip

Let $\sigma=B_1/B_2/\dots/B_k$ and $b\in B_i$. The dual of a descent
is an {\it ascent}, which is a pair $(b,B_{i-1})$ with $b>\min
B_{i-1}$. Note that this is true that each $b\in B_i$ forms an ascent because of the
canonical ordering.  So, we define the dual major index to be
$$\widehat{\mathrm{maj}}(\sigma)=\sum_{i=2}^k (i-1)(\#B_i).$$

The {\it dual inversion number} of $\sigma$, written
$\widehat{\mathrm{inv}}(\sigma)$, is the number of pairs $(b,B_j)$
such that $b\in B_i$, $b>\min B_j$, and $i>j$.  We will call these
pairs {\it dual inversions}. Clearly,
$\widehat{\mathrm{inv}}(\sigma)=\widehat{\mathrm{maj}}(\sigma)$ for
any $\sigma\in \Pi_n$, since every ascent causes $i-1$ dual
inversions.

\begin{prop} For any $\sigma\in\Pi_n$,
$$\widehat{\mathrm{inv}}(\sigma)=\widehat{\mathrm{maj}}(\sigma)=(\:\bmin{1}/\bmin{2}+\bmin{1}/\bmin{2}3)(\sigma).$$
\end{prop}

\noindent {\bf Proof:}  Let $\sigma =B_1/B_2/\dots/B_k$.  The proof
that
$\widehat{\mathrm{inv}}(\sigma)=(\bmin{1}/\bmin{2}+\bmin{1}/\bmin{2}3)(\sigma)$
is similar to the proof of Proposition 6.1. The only difference here
is that the minimum of a block can represent the $b$ in a dual
inversion $(b,B_j)$. This is taken care of by the first pattern.
$\square$

\medskip

Wachs and White \cite{WachsWhite} define four natural statistics on
partitions by encoding the partitions as restricted growth
functions.  Their statistics are {\it lb}, {\it ls}, {\it rb}, and
{\it rs}, which stand for left bigger, left smaller, right bigger
and right smaller.  For consistency, we will define these statistics
without introducing restricted growth functions, and hence the names
of the statistics may seem a little unusual.

Let $\sigma=B_1/B_2/\dots/B_k$.  If $b\in B_i$, then we will say
that $(b,B_j)$ is:

\begin{itemize}
\item a {\it left bigger pair} of $\sigma$ if $i<j$, and $b>\min B_j$,\\
\item a {\it left smaller pair} of $\sigma$ if $i>j$ and $b>\min B_j$,\\
\item a {\it right bigger pair} of $\sigma$ if $i<j$ and $b<\max B_j$,\\
\item a {\it right smaller pair} of $\sigma$ if $i>j$ and $b<\max
B_j$.
\end{itemize}

Let $lb(\sigma)$, $ls(\sigma)$, $rb(\sigma)$, and $rs(\sigma)$  be,
respectively, the number of left bigger pairs, the number of left
smaller pairs, the number of right bigger pairs, and the number of
right smaller pairs in $\sigma$.

Notice that $(b,B_j)$ is a left bigger pair if and only if it is an
inversion of $\sigma$, and $(b,B_j)$ is a left smaller pair if and
only if $(b,B_j)$ is a dual inversion of $\sigma$. Thus we have from
Propositions 6.2 and 6.3 that
$$lb(\sigma)=(\bmin{1}3/\bmin{2}\:)(\sigma),$$
$$ls(\sigma)=(\bmin{1}/\bmin{2} + \bmin{1}/\bmin{2}3)(\sigma).$$

We will now consider the other two statistics.

\begin{prop} For any $\sigma\in\Pi_n$,
$$rb(\sigma)=(\:\bmin{1}/\bmin{2}\hspace{2mm}\bmax{3}+\bmin{1}3/\bmin{2}\hspace{2mm}\bmax{4}+\bmin{1}/\bmaxmin{2}+\bmin{1}2/\bmaxmin{3})(\sigma).$$
\end{prop}

\noindent {\bf Proof:}  Let $\sigma=B_1/B_2/\dots/B_k$.  The pattern
$\bmin{1}/\bmin{2}\hspace{2mm}\bmax{3}$ counts right bigger pairs
$(b,B_j)$ where $b=\min B_i$ and $\#B_j\geq2$.  The pattern
$\bmin{1}3/\bmin{2}\hspace{2mm}\bmax{4}$ counts those pairs where
$b\not=\min B_i$ and $\#B_j\geq2$.  The other two patterns
correspond to the same two cases when $\#B_j=1$. $\square$

\medskip

The proof of the following proposition is similar to the proof of
Proposition 6.4 and is omitted.

\begin{prop}  For any $\sigma\in\Pi_n$,
$$rs(\sigma)=(\:\bmin{1}\hspace{2mm}\bmax{3}/\bmin{2}+\bmin{1}\hspace{2mm}\bmax{4}/\bmin{2}3)(\sigma).\: \square$$
\end{prop}

There has long been interest in non-crossing partitions.  Recall
that the non-crossing partitions are those in the set $\Pi_n(13/24)$
for some $n$.  Non-nesting partitions may be described as those in the
set $\Pi_n(\adj{1}{4}\hspace{-1mm}/\hspace{-1mm}\adj{2}{3})$.  Note
that this definition of a non-nesting partition is not the only one.
Klazar \cite{abbafree} defines non-nesting partitions as those in
the set $\Pi_n(14/23)$.

Recently, however, there has been increasing interest in counting
the number of crossings or nestings of a partition.  In
\cite{VasTab}, Chen et al.\ show that the crossing number and
nesting number are symmetrically distributed over $\Pi_n$ by giving
a bijection between partitions and vacillating tableaux.  In
\cite{KasZeng}, Kasraoui and Zeng give an involution of $\Pi_n$,
which exchanges the crossing number and the nesting number while
keeping another statistic, the number of alignments of two edges,
fixed.

We will describe each of these statistics and show that they too may
be translated into the language of patterns.

Let $\sigma=B_1/B_2/\dots/B_k\in\Pi_n$.  We may rewrite $\sigma$ as
a set $P\subseteq [n]\times[n]$ in the following way.  If $a,b\in
B_i$ and there is no $c\in B_i$ such that $a<c<b$ then $(a,b)\in P$.
If $B_i=\{d\}$ then $(d,d)\in P$.  It's easy to see that $P$
uniquely represents $\sigma$.  We will call $P$ the {\it standard
representation} of $\sigma$.

Let $A$ be a family $\{(i_1,j_1), (i_2,j_2)\}\subseteq P$.  We will
say that $A$ is:
\begin{itemize}
\item a {\it crossing} if $i_1<i_2<j_1<j_2$,\\
\item a {\it nesting} if $i_1<i_2<j_2<j_1$,\\
\item an {\it alignment} if $i_1<j_1\leq i_2<j_2$.
\end{itemize}

For example, the following diagram represents $\sigma=137/26/45$,
where an edge connects elements if they are adjacent in a block.

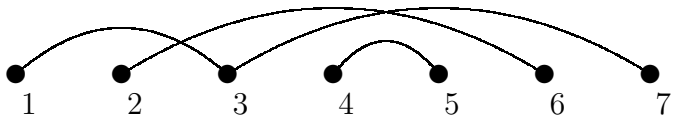
\begin{picture}(400,100)
\put(5,50){\circle*{7}} \put(45,50){\circle*{7}}
\put(85,50){\circle*{7}} \put(125,50){\circle*{7}}
\put(165,50){\circle*{7}} \put(205,50){\circle*{7}}
\put(245,50){\circle*{7}}

\put(7,35){1} \put(47,35){2} \put(87,35){3} \put(127,35){4}
\put(167,35){5} \put(207,35){6} \put(247,35){7}

\qbezier(5,50)(45,85)(85,50) \qbezier(85,50)(165,100)(245,50)
\qbezier(45,50)(125,100)(205,50) \qbezier(125,50)(145,75)(165,50)

\end{picture}

Notice that the pair $\{(1,3),(2,6)\}$ forms a crossing, the pair
$\{(2,6),(4,5)\}$ forms a nesting and the pairs $\{(1,3),(4,5)\}$
and $\{(1,3),(3,7)\}$ each form an alignment of two edges.

Let $\mathrm{cr}(\sigma)$ be the number of crossings in $\sigma$,
$\mathrm{ne}(\sigma)$ the number of nestings, and
$\mathrm{al}(\sigma)$ the number of alignments.

The following proposition is an easy consequence of the previous
definitions.

\begin{prop} For any $\sigma\in\Pi_n$,
\begin{eqnarray*}
\mathrm{cr}(\sigma)&=&(\adj{1}{3}\hspace{-1mm}/\hspace{-1mm}\adj{2}{4})(\sigma),\\
\mathrm{ne}(\sigma)&=&(\adj{1}{4}\hspace{-1mm}/\hspace{-1mm}\adj{2}{3})(\sigma),\\
\mathrm{al}(\sigma)&=&(\adj{1}{2}\hspace{-1mm}/\hspace{-1mm}\adj{3}{4}+\adj{1}{2}\adj{3}{4}+\adj{1}{2}\hspace{-1mm}\stackrel{\frown}{\hspace{1mm}{3}})(\sigma).\:
\square
\end{eqnarray*}
\end{prop}

Let $\sigma\in\Pi_n$ and $P$ be the standard representation of
$\sigma$.  Consider the family $A=\{(i_1,j_1),(i_2,j_2),\dots
(i_k,j_k)\}\subseteq P$.  Then $A$ is a $k${\it -crossing} if
$i_1<i_2<\dots<i_k<j_1<j_2<\dots<j_k$.  We say $A$ is a $k${\it
-nesting} if $i_1<i_2<\dots<i_k<j_k<j_{k-1}<\dots<j_1$.  Let
$\mathrm{cr}_k(\sigma)$ be the number of $k$-crossings of $\sigma$
and $\mathrm{ne}_k(\sigma)$ be the number of $k$-nestings of
$\sigma$. Notice that $\mathrm{cr}=\mathrm{cr}_2$ and
$\mathrm{ne}=\mathrm{ne}_2$.  The following proposition describes
these two statistics as patterns.

\begin{prop} For any $\sigma\in\Pi_n$,
\begin{eqnarray*}
cr_k(\sigma)&=&(\badj{1}{(k+1)}/\badj{2}{(k+2)}/\dots/\badji{k}{(2k)})(\sigma),\\
ne_k(\sigma)&=&(\badji{1}{(2k)}/\badj{2}{(2k-1)}/\dots/\badj{k}{(k+1)})(\sigma).\:\square
\end{eqnarray*}
\end{prop}

\section{Future Work}

There has been an explosion in interest in permutation patterns
recently, and this paper will hopefully help to generate interest in
similar work with set partitions. Sections two, three, and four
focus mainly on the question of avoidance of a partition of a three
element set, and there is more that can be done.  Klazar
\cite{abbafree, KlazarPartI, KlazarPartII}, for example, has done
work on avoidance of certain partitions of a four element set. The
problem of avoiding more than one pattern in $\Pi_4$ is yet to be
considered. Also, of interest is the problem of avoiding a family of
patterns, which include patterns from both $\Pi_4$ and $\Pi_3$.
 Sagan~\cite{Saganpartitionpatterns} has provided enumerative results for four
 different infinite families of patterns.

This is just the tip of iceberg.  We may also consider problems of
containment.  For example, what is the smallest $n$ such that we can
find a partition in $\Pi_n$, which contains all the patterns in
$\Pi_k$?  Also, for $\pi\in\Pi_k$, which $\sigma\in\Pi_n$ contain
the maximal number of copies of $\pi$?  The second question is
similar to work initiated for packing of permutations by Price in
\cite{Pricediss}.

In \cite{BabsonSteingrimsson}, Babson and Steingr{\'{i}}msson use
generalized permutation patterns to classify Mahonian statistics. It
is known that the distribution of the statistics of Milne, Sagan,
and of Wachs and White on the set of partitions of $[n]$ with $k$
blocks give nice $q$-analogues of the Stirling numbers of the second
kind. Is there any way to use the generalized patterns for set
partitions to classify the statistics which produce these nice
$q$-analogues?

Another question which arises is: what distributions do we get if we
examine these statistics on sets $\Pi_n(R)$ for some
$R\subseteq\Pi_k$?  The author is working with Sagan on a project
\cite{SagGoytqFib} that answers this question for the restricted
sets $\Pi_n(13/2)$ and $\Pi_n(13/2,123)$.  The distribution of the
statistics $ls$ and $rb$ on $\Pi_n(13/2,123)$ produce $q$-analogues
of the Fibonacci numbers, which are closely related to $q$-Fibonacci
numbers studied by Carlitz \cite{Carlitzfibnotes3, Carlitzfibnotes4} and Cigler
\cite{Ciglerqfib1}.  It is also interesting to note that these
$q$-analogues arising from restricted set partitions are related to
integer partitions.  Such $q$-analogues can also be viewed as
arising from statistics on compositions.

One partition being contained in another partition as a pattern
produces a natural partial ordering on the family of all set
partitions.  This poset is likely to be quite beautiful and have
nice structure.  It is, of course, an analogue of the poset of
permutations ordered by containment.  For more information on this
poset of permutations see \cite{Wilfpatterns}.  The author is
currently investigating properties of various posets of compositions
related to a
composition poset studied by Sagan and Vatter~\cite{SagVatcomposet}
and Bj{\"{o}}rner and Sagan~\cite{SagBjorncomposet}.

We would like to thank Ira Gessel for many informative
conversations.

\bibliography{Bib}

\end{document}